\documentclass[12pt, a4paper, twosided, reqno]{amsart}
\usepackage{enumerate, cite}
\usepackage{hyperref}
\usepackage{amsmath}
\newtheorem{theorem}{Theorem}
\newtheorem{lemma}{Lemma}

\newtheorem{example}{Example}

\newtheorem*{thA}{Theorem A}
\newtheorem*{thB}{Theorem B}
\newtheorem*{thC}{Theorem C}

\allowdisplaybreaks
\author[N. Gahlian]{Nidhi Gahlian }
\address{nidhi gahlian; department of mathematics, university of delhi, delhi-110007, india.}
\email{nidhigahlyan81@gmail.com}

\thanks {Research work of the author is supported by research fellowship from Department of Science and Technology(INSPIRE), New Delhi, India.}
\title[On Solutions ]{On Solutions of Certain Non-Linear Delay Differential Equations  }
\begin{document}
	\maketitle
	
	\begin{abstract}
		In this paper, we study the existence and non-existence of entire solutions of certain non-linear delay-differential equations.
	
	\end{abstract}
	\section{\textbf{Introduction}}
	In the study of complex differential equations, proving the existence, non existence or uniqueness of an entire or meromorphic solution for a given differential equation in the complex plane $\mathbb{C}$ is always a fascinating yet challenging task, especially when dealing with non-linear equations.\;An extension of Tumura Clunie theorem was given by Hayman\cite{hay}.\;Afterwards, exploring the characteristics of solutions of Tumura-Clunie type differential equations $f^n(z)+Q_{d}(z,f)=g(z)$ became very rigorous, many results have been proved by taking various perturbations of equations.\;One can refer to \cite{lia,liao,mues} and references therein.\;Researchers studied mainly three features i.e existence and non-existence, order of growth, and different types of forms of solutions of various types of complex  differential equations.
 
  Chen et al.\cite{chen} characterized the entire solutions of the following non-linear differential-difference equation
  \begin{equation*}
  		f^n(z)+\omega f^{n-1}(z)f'(z)+q(z)e^{Q(z)}f(z+c)=p_1e^{{\lambda}z}+e^{{-\lambda}z},
  \end{equation*}
where $n\in \mathbb{N}$,  $\omega$ and $ c \neq 0$ are complex constants.\;Moreover $q(z) \neq0, Q(z)$ are polynomials, where $Q(z)$ is not constant.\;Recently, in 2023, X. Xiang et al.\cite{xiang} investigated the growth and structure of entire solutions of the following non-linear differential-difference equation   
\begin{equation}\label{maineqA}
	f^n(z)+\omega f^{n-1}(z)f^{(k)}(z)+q(z)e^{Q(z)}f(z+c)=p_1e^{{\lambda_1}z}+e^{{\lambda_2}z},
	\end{equation}
where $n,k\geq1$ are integers, $\omega$, $p_1$, $p_2$, $\lambda_1$, $\lambda_2$, $c$ are non-zero constants, and $q\not\equiv 0$, a polynomial and $Q(z)$ is non-constant polynomial.\;The following result was proved.
 
 \begin{theorem}\cite{xiang}
 	If $f(z)$ is a transcendental entire solution with finite order to equation \eqref{maineqA}, then the following conclusions holds:
 \begin{enumerate}
 	\item If $n\geq4$ then $\rho(f)=deg(Q)=1$.
 	\item If $n\geq1$ and $\lambda(f) <\rho(f)$, then $q(z)$  degenerates into a constant, and 
 	\begin{equation*}
 		f(z)=\left(\frac{p_2 n^k}{n^k+\omega\lambda_2^k}\right)^{\frac{1}{n}}e^{\frac{\lambda_2}{n}z}, Q(z)=\left(\lambda_1-\frac{\lambda_2}{n}\right)z+\log\frac{p_1}{q\left(\frac{p_2n^k}{n^k+\omega\lambda_2^k}\right)^{\frac{1}{n}}}-\frac{\lambda_2c}{n},
 	\end{equation*}
 or
 \begin{equation*}
 		f(z)=\left(\frac{p_1n^k}{n^k+\omega\lambda_1^k}\right)^{\frac{1}{n}}e^{\frac{\lambda_1}{n}z},Q(z)=\left(\lambda_2-\frac{\lambda_1}{n}\right)z+\log\frac{p_2}{q\left(\frac{p_1n^k}{n^k+\omega\lambda_1^k}\right)^{\frac{1}{n}}}-\frac{\lambda_1c}{n}.
 \end{equation*}
 \end{enumerate}
	\end{theorem}
For the case $n=3$, they have proved the following result.
\begin{theorem}\cite{xiang}
Let $k\geq 1$ be an integer, $\omega$, $p_1$, $p_2$, $\lambda_1$, $\lambda_2$, and $c$ are non-zero constants, and $q(z)$ is a non-zero polynomial and $Q(z)$ is non-constant polynomial.\;If the non-linear delay-differential equation
\begin{equation}\label{maineqB}
f^3(z)+\omega f^{2}(z)f^{(k)}(z)+q(z)e^{Q(z)}f(z+c)=p_1e^{{\lambda_1}z}+e^{{\lambda_2}z},
\end{equation}
admits a transcendental entire solution of finite order $f$ with $\Theta(a,f)>0$,\;then $\rho(f)= deg(Q)=1$.
\end{theorem}
 In this sequence, we prove the following result for $n=3$.
 \begin{thA}\label{main1}
Suppose $\omega$, $p_1$, $p_2$, $\lambda_1$, $\lambda_2$, and $c$ are non-zero constants  with $\lambda_1 \neq\lambda_2$ and  $\frac{\lambda_1}{\lambda_2} \neq {(\frac{d}{d+1})}^{\pm 1}$ for $d=1,2$, $q(z)\not\equiv 0$, a polynomial and $Q(z)$ is non-constant polynomial.\;Then for equation \eqref{maineqB},\;there does not exist any transcendental entire solution of finite order with a finite non-zero Borel exceptional value.
\end{thA}
The following example illustrates that  it is possible for $f(z)$ to possess  zero as a Borel  exceptional value, where $f(z)$ is a transcendental entire function satisfying equation \eqref{maineqB}.
\begin{example}
 The function  $f(z)=e^{3z}$ satisfies the delay-differential equation
\begin{equation*}
	f^3(z)+3\iota f^2(z)f"(z)+4e^zf(z+\log 5)=(1+27\iota)e^{9z}+500e^{4z}.
	\end{equation*}
\end{example}
\begin{example} 
	 The function  $f(z)=e^{z}$ satisfies the delay-differential equation
\begin{equation*}
	f^3(z)+3\iota f^2(z)f"(z)+4e^zf(z+\log 5)=(1+3\iota)e^{3z}+20e^{2z}.
\end{equation*}
\end{example}
Here $f(z)$ has a zero  Borel exceptional value.\\
Now, to begin with, it is essential to note that the following outcome presented by Yang and Laine \cite{yangcc}, which stands out as a pivotal breakthrough in the context of literature on solutions for complex difference equations.
\begin{theorem}\cite{yangcc}
	Let $p(z)$, $q(z)$ be polynomials.\;Then, a non-linear difference equation
	\begin{equation*}
		f^2(z)+p(z)f(z+1)=q(z)
	\end{equation*}
	has no transcendental entire solutions of finite order.
\end{theorem}
Motivated by Theorem $3$, researchers actively engaged in extending existing theorems, and over time, the literature witnessed the emergence of more elegant results. In $2019$, Wang and Wang \cite{wang} extended and enhanced a finding by Zhang and Huang \cite{zhang}, presenting several refined outcomes within the same framework as used by Dong and Liao\cite{dong}.\;Following this, many researchers have provided varying versions of these results for certain type of non-linear differential-difference equations.
In $2022$, Wang et. al\cite{wangy} proved the slightly  changed version  of the previous result.
\begin{theorem}\cite{wangy}
	Let $n\geq3$ be an integer, $p(z)$ be a non vanishing polynomial and $p_1$, $p_2$, $\alpha_1$, $\alpha_2$ be non- zero  constants such that $\alpha_1/\alpha_2\neq n+1$ and  $\alpha_2/\alpha_1\neq n+1$.\;If the equation 
	\begin{equation}
		f^n(z)f'(z)+p(z)f(z+c)=p_1e^{\alpha_1z}+p_2e^{\alpha_2z},
	\end{equation}
admits a meromorphic solution $f(z)$ with hyper order less than $1$, then $n=3$ and $f(z)$ satisfies $\bar{\lambda}(f(z))=\rho(f(z))$.
\end{theorem} 
In this direction, a natural question that arises is what can we say when $n< 3?$  \;\; So, we replace  R.H.S with $p_1(z)e^{\lambda_1(z)}$ $+$ $p_2(z)e^{\lambda_2(z)}$ $+$ $p_3(z)e^{\lambda_3(z)}$  satisfying some condition and prove the following result.

\begin{thB}\label{theoremA}
Let $n\in\mathbb{N}$, $P(z)$ be a non vanishing polynomial. Suppose $p_j(z)$, and $\lambda_j(z)$ $(j=1,2,3)$ be non zero polynomials and non-constant polynomials.\;Also,  $\lambda_1'(z)$, $\lambda_2'(z)$, and $\lambda_3'(z)$ are distinct to each other.\;If $f$ is transcendental entire solution of the differential equation
\begin{equation}\label{maineq2}
	f^n(z)f'(z)+P(z)f(z+c)=p_1(z)e^{\lambda_1(z)}+ p_2(z)e^{\lambda_2(z)}+ p_3(z)e^{\lambda_3(z)}, 
\end{equation} 
then $f$ is of finite order but does not have zero as a Borel exceptional value.  

\end{thB}
The given example shows the necessity of $\lambda_1'(z)$, $\lambda_2'(z)$, and $\lambda_3'(z)$ should be distinct to each other.
\begin{example}
The function $f(z) = e^ z$ satisfies the differential-difference equation
\begin{equation}
	f^2(z)f'(z)+(z^2+1)f(z+\log 2)= e^{3z} + 2z^2e^z + 2e^z.
\end{equation}
\end{example}
Within the same paper, they demonstrated the following  result.
\begin{theorem}\cite{wangy}
	Let $n\geq4$ be an integer, $q(z)$ be a rational function and $Q(z)$ be a non-constant polynomial. Suppose that $c$, $p_1$, $p_2$, $\alpha_1$, $\alpha_2\in \mathbb{C}\setminus\{0\}$ with $\alpha_1\neq \alpha_2$. If the equation 
	\begin{equation}\label{6}
		f^n(z)f'(z)+q(z)f(z+c)e^{Q(z)}=p_1e^{\alpha_1z}+p_2e^{\alpha_2z},
		\end{equation}
	admits a finite order meromorphic solution $f(z)$ with $\lambda(f) <\rho(f)$ and $\lambda(\frac{1}{f}) <\rho(f)$, then $f(z)$ satsfies $\rho(f)= \deg Q =1$.\;Furthermore, $f(z)=A_1e^{\frac{\alpha_2}{n+1}z}$,\linebreak $Q(z)=(\alpha_1-\frac{\alpha_2}{n+1})z + a_0$
	or $f(z)=A_2e^{\frac{\alpha_1}{n+1}z}$, $Q(z)=(\alpha_2-\frac{\alpha_1}{n+1})z+a_0$, where $A_1,A_2, a_o \in \mathbb{C}\setminus\{0\}$.
\end{theorem}
This theorem can be extended to the case of $n=3$ by employing the same outlined techniques and fundamental principles.
The following example gives the existence for the same.\\
\begin{example}  $f(z)=e^{2z}$ is a meromorphic solution of the differential-difference equation
\begin{equation*}
	f^3(z)f'(z)+\frac{1}{8}f(z-\log 4)e^{\frac{-3}{2}z+\log16}=\frac{1}{8}e^{\frac{1}{2}z}+2e^{8z},
\end{equation*}
where $n=3$, $\alpha_1=\frac{1}{2}$, $\alpha_2=8$ and $a_0= \log16$, then $f(z)=A_1e^{\frac{\alpha_2}{n+1}z}=e^{2z}$,\linebreak $Q(z)=(\alpha_1-\frac{\alpha_2}{n+1})z + a_0= \frac{-3}{2}z+\log16$ and $\rho(f)= \deg Q=1$.
\end{example}

The subsequent example demonstrates the existence of functions for $n = 2$, although the methodology employed in Theorem $4$ is inadequate to establish it. However, there is an intuition that suggests it can be proven.
\begin{example} $f(z)=e^{z}$ is a meromorphic solution of the differential-difference equation
\begin{equation*}
	f^2(z)f'(z)+\frac{1}{7}f(z+\log 3)e^{-3z-\log3}=\frac{1}{7}e^{-z}+e^{3z},
\end{equation*}
where $n=2$, $\alpha_1=-2$, $\alpha_2=3$ and $a_0= \log3$, then $f(z)=A_1e^{\frac{\alpha_2}{n+1}z}=e^{z}$,\linebreak $Q(z)=(\alpha_1-\frac{\alpha_2}{n+1})z + a_0=-3z-\log3$ and $\rho(f)= \deg Q=1$.
\end{example}
 In $2015$, Xu et.\;al\cite{xu} considered a general differential-difference equation to obtain the following theorem.
\begin{theorem}\cite{xu}
	Consider the non-linear differential-difference equation
	\begin{equation}\label{5}
		q(z)f^n(z)+\alpha(z)f^{(k)}(z+1)=p_1(z)e^{q_1(z)}+p_2e^{q_2(z)},
		\end{equation}
	where $p_1(z),p_2(z)$ are two non zero polynomials, $q(z)$, $\alpha(z)$ are two nonzero entire functions of finite order, $q_1(z),q_2(z)$ are two non-constant polynomials and $n\geq2$ is an integer. Suppose that an entire function $f(z)$ satisfies any of the two following conditions:
		\begin{enumerate}
		\item $\lambda(f)< \rho(f)=\infty$ ; $\rho_2(f)<\infty$
		\item $\lambda_2(f)<\rho_2(f)<\infty$,
	\end{enumerate}
	then $f(z)$ can not be  a solution of equation \eqref{5}.
		\end{theorem}
Before stating  our result we first define $ P(z + w):= \sum_{i=0}^{k} b_if^{(i)}(z + w)$, where $ w \in \mathbb{C}\setminus\{0\}$ and $b_i (i = 0,1, 2, ..., k)$ are constants (atleast
one is nonzero).\;In the sequence of above, we prove the following result.
	\begin{thC}\label{theoremA}
	Let  non-linear  delay-differential equation 
	\begin{equation}\label{1}
		f^n(z)f'(z)+\alpha(z)P(z+w)=\phi(z)sin\kappa\gamma(z),
	\end{equation}
	where $n\geq2$, $\alpha(z)$ is an entire function of finite order, $\kappa$ non-zero constant, $\phi(z)$ and $\gamma(z)$ are non zero polynomial.\;If an entire function $f(z)$  satisfies any one of the following two  conditions.
	\begin{enumerate}
		\item $\lambda(f)< \rho(f)=\infty$ ; $\rho_2(f)<\infty$
		\item $\lambda_2(f)<\rho_2(f)<\infty$,
	\end{enumerate}
	then $f(z)$ can not be  a solution of equation \eqref{1}.
\end{thC}

	We must know the basic facts of Nevanlinna's value distribution theory. For a meromorphic function $f$, $n(r,f)$, $N(r,f)$, $m(r,f)$ and  $T(r,f)$ denote un-integrated counting function, integrated counting function, proximity function and characteristic function respectively. We also use first main  theorem of Nevanlinna for a meromorphic function $f$, see \cite{hay,ilpo,yanglo}.
	
	In section $2$ we state  some definitions, lemmas  and results. In section 3, we prove our main results.
	\medskip
	
	\section{\textbf{Auxiliary results}}
 We first provide elementary definitions of the order of growth $\rho(f)$,  hyper order  of growth $\rho_2(f)$, exponent of convergence of zeros $\lambda(f)$ and hyper exponent of convergence of zeros $\lambda_2(f)$ for a meromorphic function $f$.
	$$\rho(f)=\limsup_{r\to\infty}\frac{\log T(r,f)}{\log r},$$

	$$\rho_2(f) = \limsup_{r\to\infty}\frac{\log \log T(r,f)}{\log r},$$
	
	$$\lambda(f)=\limsup_{r\to\infty}\frac{\log n(r,\frac{1}{f})}{\log r}=\limsup_{r\to\infty}\frac{\log N(r,\frac{1}{f})}{\log r},$$ 
	and
		$$\lambda_2(f)=\limsup_{r\to\infty}\frac{\log \log n(r,\frac{1}{f})}{\log r}=\limsup_{r\to\infty}\frac{\log \log N(r,\frac{1}{f})}{\log r},$$ 
	where $$N(r,\frac{1}{f})=\int_{0}^{r}\frac{n(t,\frac{1}{f})-n(0,\frac{1}{f})}{t}dt + n(0,\frac{1}{f})\log r.$$ 
	A meromorphic function $g(z)$ is a small function of $f(z)$ if $T(r,g)=S(r,f)$ and vice versa.
	For a meromorphic function $f(z)$,  $S(r,f)$ denotes the quantity satisfying $S(r,f)=o(T(r,f))$, as $r\to\infty$, outside  of a possible exceptional set $E$ (not necessarily same at each occurence) of  finite linear measure. \\
	A differential polynomial $f(z)$ means it is a polynomial in $f(z)$, and its derivatives with small functions of $f(z)$ as its coefficients.\;A differential-difference polynomial of $f(z)$ means it is a polynomial in $f(z)$, its derivatives and its shifts $f(z+c)$ with small function of $f(z)$ as its coefficients.\;Borel's lemma  plays a key role in proving our all results in this paper.
	\begin{lemma}\label{imple1}(Borel's Lemma)\cite{cc}
	Suppose $f_1(z),f_2(z),...f_n(z)(n\geq2)$ are meromorphic functions and $h_1(z),h_2(z),h_3(z)\; ...\; h_n(z)$ are entire functions satisfying:
	\begin{enumerate}
		\item $\sum_{i=1}^{\infty} f_i(z)e^{h_i(z)}\equiv 0$.
		\item For $1\leq i< k\leq n$, $h_i-h_k$ are  constants.
		\item  $1\leq i\leq n$, $1\leq m< k\leq n$,  $T(r,f_i(z))=o(T(r,e^{(h_m-h_k)}))$ as $r\rightarrow \infty$ outside of a set of finite linear measure.
	\end{enumerate}
	Then $f_i\equiv0 (i=1,2,3\;...\;n)$.
	\end{lemma}
	The next lemma gives the proximity function of logarithmic derivative of a meromorphic function $f(z)$ .
\begin{lemma}\label{imple2}\cite{ilpo}
Suppose $f(z)$ is a transcendental meromorphic function and $k\geq1$ is an integer.\;Then 
\begin{equation*}
	m\left(r,\frac{f^{(k)}}{f}\right)= S(r,f).
\end{equation*}
\end{lemma}
This lemma and its various versions have made notable contributions to the progress in establishing the existence and understanding the growth characteristics of solutions in complex difference equations.\\
The next provided lemma establishes a precise asymptotic relationship between $T(r, f(z + c))$ and $T(r, f)$ for meromorphic functions of finite order only.
\begin{lemma}\label{imple3}\cite{pc}
	Suppose $f(z)$ is a meromorphic function of finie order $\rho$ and $c$ is a non-zero complex constant.\;Then for every $\epsilon>0$,
	\begin{equation*}
		T(r,f(z+c))= T(r,f)+O(r^{\rho-1+\epsilon})+O(\log r).
\end{equation*}
\end{lemma}
The next lemma gives the  discrete version of the classical logarithmic derivative of  a meromorphic function $f$ with finite order.
\begin{lemma}\label{imple4}\cite{pc}
	Suppose $f(z)$ is a meromorphic function with $\rho(f)<\infty$ and $c_1,c_2 \in \mathbb{C}$ such that $c_1\neq c_2$, then for each $ \epsilon>0$, we have
	\begin{equation*}
		m\left(r,\frac{f(z+c_1)}{f(z+c_2)}\right)=O(r^{\rho-1+\epsilon}).
	\end{equation*}
\end{lemma}
The next lemma gives the representation of entire functions having infinte order and finite hyper order.

	\begin{lemma}\label{imple5}\cite{zx}
Let $f(z)$ be a transcendental entire function of infinite order and $\rho_2(f)= \alpha < \infty $. Then $f(z)$ can be represented as
\begin{equation*}
f(z) = Q(z)e^{g(z)},
\end{equation*}
where $Q$ and $g$ are entire functions such that\\
	$\lambda(Q) = \rho(Q) = \lambda(f),$ $ \lambda	_2(Q) = \rho_2(Q) = \lambda_2(f) $ and		$\rho_2(f) = \max \{\rho_2(Q), \rho_2(e^g)\}$.
	\end{lemma}
	\section{\textbf{Proof of main theorem}}
\begin{proof}[\textbf{\underline{Proof of Theorem A:}}]
We will establish the outcome  by contradiction. Suppose $f(z)$ is finite order transcendental entire function satisfying equation\eqref{maineqB}  with a non-zero finite Borel exceptional value. We proceed to investigate the three cases outlined below.\\
\textbf{Case 1:} If $\rho{(f)}>1$ and let $\alpha$ be a non-zero finite Borel exceptional value. Then we can utilize Weierstrass factorization theorem to represent $f(z)$ as 
\begin{equation}\label{eqc}
f(z)=r(z)e^{s(z)} +\alpha,
\end{equation}
where $r(z)$ is an entire function such that $\rho(r)<$ deg $s(z)$ and $s(z)= a_tz^t+a_{t-1}z^{t-1}+...a_0, a_t\neq0$ is a polynomial of degree $t>1$.
From equation \eqref{maineqB}, we have 
\begin{align*}
    m(r,e^{Q(z)})&= m\left(r,\frac{q(z)e^{Q(z)}f(z+c)}{q(z)f(z+c)}\right)\\&=
m\left(r,\frac{p_1e^{\lambda_1z}+p_2e^{\lambda_2z}-f^2(z)(f(z)+\omega f^{(k)}(z))}{q(z)f(z+c)}\right)\\&\leq m\left(r,\frac{p_1e^{\lambda_1z}+p_2e^{\lambda_2z}}{q(z)f(z+c)}\right)+ m\left(r,\frac{f^2(z)(f(z)+\omega f^{(k)}(z))}{q(z)f(z+c)}\right) + O(1)\\& \leq
m(r,p_1e^{\lambda_1z}+p_2e^{\lambda_2z})+m\left(r,\frac{1}{q(z)f(z+c)}\right)+m(r,f^2(z))+\\& \ \ \ m\left(r,\frac{f(z)}{q(z)f(z+c)}\right)+m\left(r,\frac{f^{(k)}(z)}{q(z)f(z+c)}\right) +O(1)
	\end{align*}
Applying lemma \eqref{imple2} and \eqref{imple4} in the above inequality, we get
\begin{align*}
 T(r,e^{Q(z)}) &\leq m\left(r,\frac{1}{q(z)f(z+c)}\right)+2T(r,f(z))+ m\left(r,\frac{f^{(k)}(z)}{q(z)f(z+c)}\right)+ S(r,f)\\& 
\leq m\left(r,\frac{1 f(z)}{q(z)f(z)f(z+c)}\right)+2T(r,f(z))+ m\left(r,\frac{f^{(k)}(z)f(z)}{q(z)f(z+c)f(z)}\right)\\&\leq 3T(r,f(z))+S(r,f),
\end{align*}
 with the help of simple calculations. Above inequality implies $ deg (Q(z))\leq \rho(f)=t$. Also, $ deg(Q(z))\geq 1$  by hypothesis. Next we discuss it in two parts;

 	\textbf{Case 1(i):} Suppose $1\leq deg(Q(z))<t$, then  replacing the function $f(z)$  in  equation \eqref{maineqB}  with its corresponding expression from equation \eqref{eqc}.
 	we get, 
 	\begin{align*}
  p_1e^{\lambda_1z}+p_2e^{\lambda_2z}=&r^3(z)e^{3s(z)} +3\alpha^2r(z)e^{s(z)}+ 3\alpha r^2(z)e^{2s(z)}+ \omega \alpha^2 e^{s(z)}\kappa(z)+\\&  2 \omega \alpha r(z)e^{s(z)}e^{s(z)}\kappa(z) +\omega r^2(z)e^{2s(z)}e^{s(z)}\kappa(z)+ \alpha q(z)e^{Q(z)}\\& \  \ +q(z)e^{Q(z)}r_c(z)e^{s_c(z)} +\alpha^3,
 	\end{align*}
 	where $r_c(z)= r(z+c), s_c(z)=s(z+c),e^{s(z)}\kappa(z)=f^{(k)}$ and $\kappa(z)$ is as follows:
 	\begin{align*}
 	f(z)&=r(z)e^{s(z)} +\alpha,\\
 	f'(z)&=r'(z)e^{s(z)} +r(z)s'(z)e^{s(z)}+0\\&=e^{s(z)}[r'(z)+r(z)s'(z)].
 	\end{align*}
 Here $f'(z)$ is the product of the function $e^{s(z)}$ and a small function of $e^{s(z)}$.  Continuing in this way, we get $k^{th}$ derivative of $f(z)$ is the multiplication of the function $e^{s(z)}$ and some small function of $e^{s(z)}$ say $\kappa(z)$.
 Now, from above equation we get
 \begin{align*}
   p_1e^{\lambda_1z}+p_2e^{\lambda_2z}=&e^{3s(z)}(r^3(z) +\omega r^2(z)\kappa (z)) + e^{2s(z)}(3r^2\alpha +2\omega \alpha r(z)\kappa(z))\\&  +e^{s(z)}(3\alpha^2r(z)+\omega \alpha^2 \kappa(z))+\\& \alpha q(z)e^{Q(z)}+  q(z)e^{Q(z)}r_c(z)e^{s_c(z)}.
\end{align*}
Taking, 
\begin{align*}
	&	D_1(z)=r^3(z) +\omega r^2(z)\kappa (z),\\&
	D_2(z)=3r^2\alpha +2\omega \alpha r(z)\kappa(z),\\&
	D_3(z)=3\alpha^2r(z)+\omega \alpha^2 \kappa(z),
\end{align*}
 and simplifying we get,
\begin{align}\label{eqd}
D_1(z)e^{3s(z)}+D_2(z)e^{2s(z)}+(D_3(z)e^{s_1(z)} +q(z)r_c(z)e^{\beta(z)})&e^{a_tz^t}+\alpha^3+\alpha q(z)e^{Q(z)}-\nonumber\\&  p_1e^{\lambda_1z}-p_2e^{\lambda_2z}=0,
	\end{align}
where $s_1(z)=a_{t-1}z^{t-1}+a_{t-2}z^{t-2}...+a_0$, $\beta(z)=Q(z) +a_t(t_{c_1}z^{t-1}+t_{c_2}z^{t-2}+ ...c^t)+ a_{t-1}(z+c)^{(t-1)}+...a_0$ are polynomials of degree atmost $t-1$.
Applying lemma \ref{imple1}  in equation \eqref{eqd}, we get
\begin{equation}\label{eqe}
	\alpha q(z)e^{Q(z)}+ \alpha^3-p_1e^{\lambda_1z}-p_2e^{\lambda_2z}=0.
\end{equation}
Again apply lemma \ref{imple1} on equation \eqref{eqe}. We have three cases.

\begin{enumerate}[(i)]
	\item If $Q(z)\neq\lambda_1\neq \lambda_2$, then we have $p_1\equiv0 \equiv \alpha^3\equiv p_2$, which is contradiction to the given hypothesis.
	
\item  If $Q(z)=\lambda_2$, then we have $p_1\equiv0 \equiv \alpha^3$, which is again contradiction to the given hypothesis.
	
   \item If $Q(z)=\lambda_1$, then we have $p_2\equiv0 \equiv \alpha^3$, which is also a contradiction to the given hypothesis.
   \end{enumerate}
so  $1\leq deg(Q(z))<t$ is not possible.\\
\textbf{Case 1(ii):} If $deg(Q(z))=t$ say $Q(z)=\gamma_tz^t+\gamma_{t-1}z^{(t-1)}+...\gamma_0$, $\gamma_t\neq0$ is a polynomial and
equation,
\begin{align*}
D_1(z)e^{3s(z)}+D_2(z)e^{2s(z)}+  D_3(z)e^{s(z)} +q(z)r_c(z)e^{Q(z)}&e^{s(z+c)}+\alpha q(z)e^{Q(z)} +\alpha ^3-\\& p_1e^{\lambda_1z}-p_2e^{\lambda_2z}=0,
\end{align*}
implies 
\begin{align}\label{eqf}
   D_1(z)e^{3s_1(z)}e^{3a_tz^t}+ &D_2(z)e^{2s_1(z)}e^{2a_tz^t} +q(z)r_c(z)e^{Q_1(z)+s_1(z+c)}e^{\gamma_tz^t +a_t(z+c)^t}\nonumber\\& + D_3(z)e^{s_1(z)}e^{a_tz^t}+ \alpha q(z)e^{Q_1(z)}e^{\gamma_tz^t} +\alpha^3-p_1e^{\lambda_1z}-p_2e^{\lambda_2z}=0,
\end{align}
where  $Q_1(z)=\gamma_{t-1}z^{t-1}+...\gamma_0$ and  $s_1(z)$ is same as defined  above. Now discuss the following subcases:\\
\textbf{ 1(ii)(a):}If $\gamma_t\neq Ka_t$ for $K=1,2,3.$ Apply Borel's lemma \ref{imple1} on equation \eqref{eqf}, we have $\alpha q(z)e^{Q_1(z)}\equiv0$ implying either $\alpha\equiv0$ or $q(z)\equiv0$, which is not possible.\\
\textbf{ 1(ii)(b):}If $\gamma_t= Ka_t$ for some $K=1,2,3,$ say $\gamma_t=2a_t$, then equation \eqref{eqf} becomes
\begin{align*}
(D_1(z)e^{3s_1(z)}+q(z)r_c(z)e^{\beta (z)})e^{3a_tz^t}+ &(D_2(z)e^{2s_1(z)} +\alpha q(z)e^{Q_1(z)})e^{2a_tz^t}+D_3e^{s_1(z)}e^{a_tz^t}\\&- -p_1e^{\lambda_1z}-p_2e^{\lambda_2z}+\alpha^3=0,
\end{align*}
where $\beta(z)=Q_1(z)+s_1(z+c)+a_t[\binom t1 cz^{t-1}+\binom t2 c^2z^{t-2}+...c^t]$ is a polynomial of degree  atmost $t-1$. Now applying Borel's lemma \ref{imple1}, we have
$p_1e^{\lambda_1}+p_2e^{\lambda_2 z}-\alpha^3=0$, which is also not possible.\\
\textbf{1(ii)(c):} If $\gamma_t=-Ka_t$ for some $K=1,2,3,$ say $\gamma_t=-a_t$, then equation \eqref{eqf} becomes
\begin{align*}
D_1(z)e^{3s_1(z)}e^{3a_tz^t}+ D_2(z)e^{2s_1(z)}e^{2a_tz^t}  + &D_3(z)e^{s_1(z)}e^{a_tz^t} + q(z)r_c(z)e^{\beta(z)} +\alpha^3-p_1e^{\lambda_1z} \\& +\alpha q(z)e^{Q_1(z)}e^{-a_tz^t} -p_2e^{\lambda_2z}=0,
\end{align*}
where $\beta(z)$ is same as defined above. Applying  Borel's lemma \ref{imple1} on above equation, we get $\alpha q(z)e^{Q_1(z)}\equiv0$ which further implies  either $\alpha=0$ or $ q(z)\equiv0$ and it is  contradictory to our hypothesis.

\textbf{Case 2:} If $\rho(f)<1$ then, using similar technique as in [\cite{xiang},Theorem 1.1],we get a contradiction.

\textbf{Case 3:} If $\rho(f)=1$ then, represent   $f(z)=U(z)e^{Vz+g}+\alpha$, where $V\neq0,g \in \mathbb{C} $ $\&$ $ U(z) $ is an entire function such that $\rho(U)<1$. Proceding with the same steps as in case $1$, we have $deg(Q(z))=1$.
So $ Q(z)=hz+p$, where
 $h\neq0$ $\&$ $p \in \mathbb{C}$. Substituting values of $Q(z)$ and $f(z)$ in equation \eqref{maineqB}, we have
\begin{equation*}
	\begin{split}
 &U^3(z)e^{3(Vz+g)}+\alpha^3+ 3\alpha^2 U(z)e^{Vz+g}+3\alpha U^2(z)e^{2(Vz+g)}+\alpha q(z)e^{hz+p} -p_1e^{\lambda_1z}-p_2e^{\lambda_2z}\nonumber \\&\omega ( U^2(z)e^{2(Vz+g)}+\alpha^2+2\alpha U(z)e^{Vz+g})e^{Vz+g}\kappa'(z) + q(z)e^{hz+p}U_c(z)e^{Vz+Vc+g} =0,
 \end{split}
\end{equation*}
where $ U(z+c)=U_c(z) $ and $e^{Vz+g}\kappa'(z)=f^{(k)}$, $\kappa'(z)$ is small function of  $e^{Vz+g}$.
Above equation can be written as 
\begin{equation}\label{eqg}
	\begin{split}
G_1e^{3Vz}+ G_2e^{2Vz}+G_3e^{Vz}+G_4e^{(h+V)z}+&\alpha q(z)e^pe^{hz}+\alpha^3- p_1e^{\lambda_1z}-p_2e^{\lambda_2z}=0,
\end{split}
\end{equation}
where \begin{align*}
	&	G_1(z)=(U^3(z) +\omega U^2(z)\kappa' (z))e^{3g},\\&
	G_2(z)=(3U^2(z)\alpha +2\omega \alpha U(z)\kappa'(z))e^{2g},\\&
	G_3(z)=(3\alpha^2U(z)+\omega \alpha^2 \kappa'(z))e^g,\\&
	G_4(z)=q(z)U_c(z)e^{vc+g+p}.
\end{align*}
Here also we discuss the following subcases.\\
\begin{enumerate}[(i)]
\item  If $KV\neq h$ and $KV\neq \lambda_1\neq \lambda_2$ for any $K=1,2,3$. Then applying lemma \ref{imple1} on equation \eqref{eqg}, we get $p_1 \equiv 0\equiv p_2$, which  is not possible.
\item If $KV= h$ for some $K=1,2,3$ but $KV\neq \lambda_1\neq \lambda_2$ for any $K=1,2,3$. Then applying lemma \ref{imple1} on equation \eqref{eqg}, we obtain  $p_1 \equiv 0\equiv p_2$, which  is again not possible.
\item If $KV\neq h$, $KV\neq \lambda_1$ for any $K=1,2,3$,  but $KV=\lambda_2$ for some $K=1,2,3$. Then applying lemma \ref{imple1} on equation \eqref{eqg}, we obtain  $p_1 \equiv 0$, which is a contradiction .
\item If $KV\neq h$ for any $K=1,2,3$,  but  $KV=\lambda_1$ and $KV=\lambda_2$ for some $K=1,2,3$. Then it will not be the case as it violates the assumptions of $\lambda_1$ and $\lambda_2$. So, $KV=\lambda_1$ and $KV=\lambda_2$ for some $K=1,2,3$ will never be a case.
\item  If $KV\neq h$, $KV\neq \lambda_2$ for any $K=1,2,3$,  but $KV=\lambda_1$ for some $K=1,2,3$. Then applying lemma \ref{imple1} on equation \eqref{eqg}, we obtain  $p_2 \equiv 0$, which is a contradiction .
\item If $KV=h$ and $KV=\lambda_1$  for some $K=1,2,3$,  but $KV\neq \lambda_2$ for any $K=1,2,3$.Then applying lemma \ref{imple1} on equation \eqref{eqg}, we obtain  $p_2 \equiv 0$, which is also a contradiction. 
\item  If $KV=h$ and $KV=\lambda_2$  for some $K=1,2,3$ ,  but $KV\neq \lambda_1$ for any $K=1,2,3$. Then applying lemma \ref{imple1} on equation \eqref{eqg}, we obtain  $p_1 \equiv 0$, which is also a contradiction. \\
Hence, the proof.
\end{enumerate}
\end{proof}
\begin{proof} [\textbf{\underline{Proof of Theorem B:}}]
	Let $f$ be  a transcendental entire solution that satisfies equation \eqref{maineq2}, then from equation \eqref{maineq2} and lemma \ref{imple3} we have,
		\begin{align*}
	f^n(z)f'(z)&=p_1(z)e^{\lambda_1(z)}+ p_2(z)e^{\lambda_2(z)}+ p_3(z)e^{\lambda_3(z)}-P(z)f(z+c)&\\
(n+1)T(r,f)&\leq T(r,p_1(z)e^{\lambda_1(z)}+ p_2(z)e^{\lambda_2(z)}+ p_3(z)e^{\lambda_3(z)}-P(z)f(z+c))&\\
&\leq T(r,p_1(z)e^{\lambda_1(z)}+ p_2(z)e^{\lambda_2(z)}+p_3(z)e^{\lambda_3(z)})+ T(r,P(z))+\\& \ \ \  \ \ \ T(r,f(z+c))+log3.
\end{align*}
	Using first fundamental theorem of Nevanlinna theory, we get 
	\begin{equation}\label{eqa}
		nT(r,f)\leq Sr^\gamma+ S(r,f),
		\end{equation}
	where $\gamma:= \max \{\deg\lambda_1(z), \deg\lambda_2(z), \deg\lambda_3(z)\}$\\
	 and \\
	$$ S:=\frac{\mbox{sum of leading coefficients of}\  \lambda_1(z),\lambda_2(z),\lambda_3(z)}{\pi}.$$
Hence,  equation \eqref{eqa} implies $f$ has a finite order.
\smallskip 
	 
	Now, we prove the remaining  part of the  theorem by contradiction. Assume  $f(z)=\alpha(z)e^{\beta(z)}$ be a solution of equation \eqref{maineq2}, where $\alpha(z)$ ia a non zero polynomial and $\beta(z)$ is a non constant polynomial. Setting $f(z)$ in  equation \eqref{maineq2} gives,
	\begin{align*}
		p_1(z)e^{\lambda_1(z)}+p_2(z)e^{\lambda_2(z)}+p_3(z)e^{\lambda_3(z)}&=f^n(z)f'(z)+P(z)f(z+c)
		\\&=(\alpha^n(z)e^{n\beta(z)})(\alpha'(z)e^{\beta(z)}+\beta'(z)e^{\beta(z)}\alpha(z))+\\
		& \  \  \  \ P(z)\alpha(z+c)e^{\beta(z+c)}.
	\end{align*}
	We have,
	\begin{equation}\label{maineq3}
		p_1(z)e^{\lambda_1(z)}+p_2(z)e^{\lambda_2(z)}+p_3(z)e^{\lambda_3(z)} -R(z)e^{(n+1)\beta(z)}-P(z)\alpha_c(z)e^{\beta_c(z)}=0,
	\end{equation}
	where  $ R(z)=\alpha^n(z)\alpha'(z)+\beta'(z)\alpha^{(n+1)z}$, $\alpha_c(z)=\alpha(z+c)$ and $\beta_c(z)=\beta(z+c)$.\\
	Now we discuss the possible situations, as $\beta(z)-(n+1)\beta(z)$ is never constant, so we have only four possible cases.\\
	\begin{enumerate}
 \item If $\beta_c(z)-\lambda_j(z) \neq$ constant and $(n+1)\beta(z)-\lambda_j(z)\neq$ constant for any $j=1,2,3.$ Now applying Borel's lemma \ref{imple1} in equation \eqref{maineq3}, we have
 $p_1(z)\equiv0\equiv p_2(z)\equiv p_3(z) $ which is contradictory to our assumption.\\
 \item If  $\beta(z)-\lambda_j(z)\neq$ constant for any $j=1,2,3,$ and      $(n+1)\beta(z)-\lambda_j(z)=$ constant, for some $j=1,2,3$, say  $(n+1)\beta(z)-\lambda_1(z)=$ contant. Let
 \begin{equation*}
 (n+1)\beta(z)= \alpha_sz^s+\alpha_{s-1}z^{s-1}+...\alpha_0
 \end{equation*}
and
\begin{equation*}
	\lambda_1(z)=\alpha_sz^s+\alpha_{s-1}z^{s-1}+...\lambda_0,
 \end{equation*}
such that $(n+1)\beta(z)-\lambda_1(z)=\alpha_0-\lambda_0=$ constant. Then, equation \eqref{maineq3} becomes
\begin{equation*}
 \; \; \; \; \; \;  (p_1(z)e^{\lambda_0}-R(z)e^{\alpha_0})e^{\gamma(z)}+p_2(z)e^{\lambda_2(z)}+p_3(z)e^{\lambda_3(z)}-P(z)\alpha_c(z)e^{\beta_c(z)}=0,
\end{equation*} 
where $\gamma(z)=\alpha_sz^s+\alpha_{s-1}z^{s-1}+...\alpha_1z$. Applying Borel's lemma\ref{imple1}, we have
$p_2(z)\equiv0$ and $p_3(z)\equiv0$, which is again a contradiction.

\item If $\beta(z)-\lambda_j(z)=$ constant for some $j=1,2,3,$ say  $\beta(z)-\lambda_1(z)=$ constant, and  $(n+1)\beta(z)-\lambda_j(z)\neq$ constant for any $j=1,2,3$. Then by doing same calculations as in previous case and applying Borel's lemma \ref{imple1}, we get the same contradiction.
 \item If  If $\beta(z)-\lambda_j(z)=$ constant for some $j=1,2,3,$ say  $\beta(z)-\lambda_1=$ constant, and  $(n+1)\beta(z)-\lambda_j(z)=$ constant for some $j=1,2,3$, say $(n+1)\beta(z)-\lambda_2(z)=$ constant. Applying Borel's lemma\ref{imple1}. into equation \eqref{maineq3}, we have
  $p_3(z)\equiv0$, which is also not possible.
\end{enumerate}	
\end{proof}
\begin{proof}[\textbf{\underline{Proof of Theorem C:}}]
	
(1) Let $f(z)$ be an entire solution of equation \eqref{1} satisfying $\lambda(f)< \rho(f)=\infty$ ; $\rho_2(f)<\infty$ . By lemma \ref{imple5}, we can write
	\begin{equation*} 
		f(z)=U(z)e^{V(z)},
	\end{equation*}
	where $U(z)$ is an entire fuction, and $V(z)$ is a transcendental entire function such that, \begin{equation}\label{A}
		\lambda(U)=\rho(U)=\lambda(f), 
		\lambda_2(U)=\rho_2(U)=\lambda_2(f).
	\end{equation} Also,
	\begin{equation*}
		\rho_2(f)= \max \{ \rho_2(U),\rho_2(e^V)\}.
	\end{equation*} 
	Now, with the given assumption $\rho_2(f)<\infty$, we have
	\begin{equation*}
		\rho(V(z))=\rho_2(e^{V(z)})<\infty.
	\end{equation*}  
	Substituting $f(z)=U(z)e^{V(z)}$ in equation  \eqref{1},we have
	\begin{align}\label{2}
		\phi(z)\sin\kappa\gamma(z) =&\; U^n(z)e^{nV(z)}(U'(z)e^{V(z)}+ U(z)V'(z)e^{V(z)}) \nonumber \\ &  +\alpha(z)\sum_{j=1}^{k}(U(z+w)e^{V(z+w)})^{(j)}\nonumber\\
		 \phi(z)\sin\kappa\gamma(z)=&\; e^{(n+1)V(z)}U^n(z)(U'(z)+U(z)V'(z))+\alpha(z)S(z)e^{V(z+w)},
	\end{align}
	where $S(z)$ is a differential polynomial of $U(z+w)$ and $ V(z+w$) with $\rho(S)<\infty$.
	Let $M(z)=V(z+w)-(n+1)V(z)$, then above equation becomes
	\begin{align*}
		e^{(n+1)V(z)}U^n(z)(U'(z)+U(z)V'(z))+\alpha(z)S(z)e^{M(z)+(n+1)V(z)}	=\phi(z)\sin\kappa\gamma(z).
	\end{align*}
	Dividing above equation  by $e^{(n+1)V(z)}$, we get
	\begin{equation*}
		U^n(z)(U'(z)+U(z)V'(z))+\alpha(z)S(z)e^{M(z)}=\phi(z)\sin\kappa\gamma(z)e^{-(n+1)V(z)}.
	\end{equation*}
	Now, there arises two cases.
	\begin{enumerate}[(i)]
		\item If $ M(z)$ is a polynomial.
		\begin{equation*}
			\rho(U^n(z)(U'(z)+U(z)V'(z))+\alpha(z)S(z)e^{M(z)})<\infty,
		\end{equation*}
		and hence,
		\begin{equation*}
			\rho(\phi(z)\sin\kappa\gamma(z)e^{-(n+1)V(z)})<\infty,
		\end{equation*}
		but it contradicts with the hypothesis implying $M(z)$ can't be a polynomial.
	\item  If $ M(z)$ is a transcendental entire function, then equation \eqref{2} can be written as 
		\begin{equation}
			\begin{split}	&e^{(n+1)V(z)}U^n(z)(U'(z)+U(z)V(z))+\alpha(z)S(z)e^{V(z+w)}\\&-e^{\beta(z)}\phi(z)\sin\kappa\gamma(z)=0,
				\end{split}
			\end{equation}
			where $\beta(z) \equiv 0$. Utilizing lemma \ref{imple1}, we deduce that 
			\begin{align*}
				U^n(z)(U'(z)+U(z)V(z))\equiv0, \alpha(z)S(z)\equiv0, \phi(z)\sin\kappa\gamma(z)\equiv0.
			\end{align*}
			If $U^n(z)(U'(z)+U(z)V(z))\equiv0$ , $U(z)\not\equiv0$. On solving it, we obtain
			\begin{equation}
				U(z)=ce^{-V(z)}.
			\end{equation} 
		This implies $\rho(U(z))=\infty$, which contradicts with equation \eqref{A} and the given assupmtion.
		\end{enumerate}
		Hence $f(z)$ with $\lambda(f)< \rho(f)=\infty$ ; $\rho_2(f)<\infty$ is not the solution of equation \eqref{1}.\\ \\
		Again, let  $f(z)$ be an entire solution of equation \eqref{1} with $\lambda_2(f)<\rho_2(f)<\infty$. By lemma \ref{imple5}, we can write
		\begin{equation*} 
			f(z)=R(z)e^{Q(z)},
		\end{equation*}
		where $R(z)$ is an entire fuction, and $Q(z)$ is a transcendental entire function such that, \begin{equation}\label{A}
			\lambda(R)=\rho(R)=\lambda(f), 
			\lambda_2(R)=\rho_2(R)=\lambda_2(f).
		\end{equation} Also, 
		\begin{equation*}
			\rho_2(f)= \max \{ \rho_2(R),\rho_2(e^Q)\}.
		\end{equation*} 
		According to the given conditions
		\begin{equation*}
			\rho_2(f(z))=\rho_2(e^{Q(z)})<\infty,
		\end{equation*} 
		so $\rho_2(R(z))<\rho_2(e^{Q(z)})=\rho(Q(z))<\infty$ and hence $\rho_2(Q(z))=0$. Now substituting $f(z)=R(z)e^{Q(z)}$ in equation \eqref{1}, and proceeding with same steps as in above, we obtain
		\begin{equation}
			R^n(z)(R'(z)+R(z)Q'(z))+\alpha(z)S(z)e^{N(z)}=\phi(z)\sin\kappa\gamma(z)e^{-(n+1)Q(z)},
		\end{equation}
		where $N(z)=Q(z+w)-(n+1)Q(z)$.
		\begin{enumerate}[(i)]
			\item If $\rho(N(z))<\rho(Q(z)$, then
			\begin{align*}
\rho_2(R^n(z)(R'(z)+R(z)Q'(z))+\alpha(z)S(z)e^{N(z)})&\leq\max\{\rho_2(R(z)),\rho(N(z))\}\\&<\rho(Q(z))\\&=\rho_2(\phi(z)\sin\kappa\gamma(z)e^{-(n+1)Q(z)}),
			\end{align*}
	which is a contradiction.		
			\item If $\rho(N(z))=\rho(Q(z)$, applying the same steps as in part $(1)$ and using lemma \ref{imple1} we get a contradiction.
		\end{enumerate}
	Hence, $f(z)$ with $\lambda_2(f)<\rho_2(f)<\infty$ is not the solution of equation \eqref{1}.
	\end{proof}

\section*{\textbf{Acknowledgement}}
I am thankful to my supervisor Professor Sanjay Pant for his valuable comments and suggestions to improve the readibility of the paper.

\end{document}